\documentclass[12pt]{amsart}
\usepackage{type1ec}
\usepackage{latexsym}
\usepackage{amsmath}
\usepackage{amssymb}
\usepackage{mathrsfs}
\usepackage{graphicx}
\input epsf

\newcommand{\dom}{\operatorname{dom}}
\newcommand{\im}{\operatorname{im}}
\newcommand{\On}{\operatorname{On}}

\begin{document}

\title
[Nonstandard Analysis: Its Creator and Place]
{
NONSTANDARD ANALYSIS:\\ ITS CREATOR AND PLACE
}
\dedicatory{On the 95th anniversary of the birth of Abraham Robinson}
\author{S.~S. Kutateladze}
\date{June 15, 2013}
\begin{abstract}
This is a biographical sketch and  tribute to Abraham Robinson
on the 95th anniversary of his birth with a short discussion of the place of nonstandard
analysis in the present-day mathematics.
\end{abstract}
\address[]{
Sobolev Institute of Mathematics\newline
%\indent 4 Koptyug Av.\newline
\indent Novosibirsk 630090, RUSSIA
}

\email{
sskut@math.nsc.ru
}
\maketitle

This year, the world mathematical community recalls the memory of Abraham Robinson 
(1918--1974),
an outstanding scientist whose contributions to delta-wing theory and model theory
are the most convincing proofs of unity between pure and applied mathematics.
Robinson created nonstandard analysis which is one of the most controversial, marvelous,
and intriguing applications of logic to the core of mathematics.

\section*{The Life's Signposts of Abraham Robinson}

Abraham Robinson was born on October 6, 1918  in  Lower Silesia
at a small Prussian town Waldenburg
(today this is  Wa{\l}brzych in Poland).\footnote{Cp.~\cite{Mostow}--\cite{Dauben} for biographical details.}
In America Abraham was lately abbreviated as Abby.

Abraham, Abby,  received the name in honor of his father who had died
young not long before the birth of his younger son.
The surname was written ``Robinsohn'' those days.
Abby's father was a~hebraist,\footnote{A~hebraist is an expert in Hebrew Studies.} talmudist, and zionist. Abby's grandfather on the mother's
side was also a~talmudist.
Abby's uncle Isac was a famous and successful surgeon, and Abby together with his elder
brother Saul spent summers at Isac's home near Vienna.

In~1925  Abby's mother, Lotte Robinsohn, moved with her two boys to
Breslau, the capital of Silesia, where there was a large Jewish community.
The brothers learned in a  Jewish private school founded by
Rabbi Max Simonson, who took great care of the junior Robinsohns, remarking
that ``the big boy was an extremely gifted child, but the little one was a genius.''

 After Hitler seized power in~1933, Robinsohns emigrated to Palestine.
The family  settled in Jerusalem, where Abby went to the Rehavia Secondary School,
appraised with the excellency of his Hebrew. He and his brother enlisted into  Haganah,
the illegal organization for defence against Arabs. In due time Haganah turned into
a basis of the Army of Defense for Israel (Tzahal).

In~1936 Robinson entered the Einstein Mathematical Institute which was
actually the Mathematics Department of the Hebrew University of Jerusalem.
His mentor was Abraham Fraenkel.

In 1938 the first paper of Abby appeared in the {\it Journal of  Symbolic Logic}
(with h in the author's surname). In January of 1940  Robinson and his companion
Jacob Fleischer moved to Paris where Abby enrolled  in~Sorbonne.
But the Germans occupied Paris in June, and Abby together with Fleischer
escaped to England through Bordeaux. In England he joined
the Free French that was collected by de Gaulle. He became a sergeant of the
Free French Air Force. Although Abby was a subject of the British
Crown his German origins had hampered his entrance into the British Army for a time being.
Miraculously enough, Abby's merits came into play. He had helped a familiar French Captain
to make  a memorandum on aircraft wings for the Ministry of Aircraft Production,
and he was soon reassigned to the British Air Force and transferred  to the Royal
Aircraft Establishment in Farnborough as an assistant (grade~3) in the British
Ministry of Aircraft Production.

In December of 1942  Robinson wrote to his supervisors in Jerusalem that he
had decided to participate in the general struggle against Fascism and apply his knowledge in
applied mathematics to this end. He remarked that there was no effort for him
to turn to applied problems. Robinson addressed  the problem of comparison between
single-engine and twin-engine planes for which he  suggested an analog of the
variational method by Ludwig Prandl. He also worked on the problem of structural fatigue and collapse of a flying boat.

In~1944 Robinson married Ren\'ee Kopel, a fashion photographer.
Abby lived with Ren\'ee up to his terminal day.

Robinson was a member of the group studying the German V-2 missiles
as well as a mission of the British  Intelligence Objectives Subcommission
which concerned intelligence gathering about the aerodynamical research in Germany.
In~1946 Robinson returned to Jerusalem to pass  examinations
for the Master degree. The results were as follows: ``physics good,
mathematics excellent.'' In this short period  Abby worked together
with Theodore Motzkin.

In~1946 the Royal College of Aeronautics was founded in Cranfield near London.
Robinson was offered the position of a Senior Lecturer with salary 700 pounds per year.
It is worth mentioning that Robinson was the only member of the teaching staff
who learned how to pilot a plane.
In Cranfield Abby became a coauthor of delta-wing theory for supersonic flights, and
in 1947 he learned Russian in order to read the Soviet scientific periodicals.

To gain the  PhD degree, Robinson joined the Birkbeck College which was intended for
mature working students and provided instructions mainly in the evening or on weekends.
Abby's supervisor in the college was Paul Dienes, a Hungarian specialized mainly in
function theory.  Dienes instigated Abby's interest in summation methods (which
resulted lately in Abby's work with Richard Cooke who also taught in the Birkbeck College).
Dienes was a broad-minded scientist with interests in algebra and foundations. In 1938
he published the book {\it Logic of Algebra}, the topic  was close to Abby's train
of thoughts. In this background Robinson returned to logic and presented and maintained
the PhD thesis ``On the Metamathematics of Algebra'' in~1947.

   In~1951 Robinson moved to Canada where he worked at the
Depart\-ment of Applied Mathematics of Toronto University.
He delivered lectures on differential equations, fluid mechanics, and
aerodynamics. He also supervises postgraduate students in applied  mathematics.
Abby worked on similarity analysis and wrote
``Foundations of Dimensional Analysis'' which was
published only after his death in~1974.

Robinson was the theorist of delta-wing, but his Farnborough research in the area
was highly classified. In Toronto Robinson
wrote his {\it magna opus} in aerodynamics, {\it Wing Theory}, which was based
on the courses he delivered in Cranfield as well as on his research in~Canada.
Robinson invited  as a coauthor John Laurmann, his former student in Cranfield.
The book addressed airfoil design of wings under subsonic and supersonic speeds
in steady and unsteady flow conditions.
James Lighthill, the creator of aeroacoustics and  one of the most prominent mechanists of the twentieth century, appraised most of the book as ``an admirable compendium of the mathematical theories of the aerodynamics of air\-foils and wings.''
Robinson performed some impressive studies of aircraft icing and waves in elastic media,
but in the mid-1950s his interest in applied topics diminished had been fading.
Robinson continued lecturing on applied mathematics, but arranged a seminar of logic for a small group of curious  students.

In~1952 Robinson participated in the second Colloquium on Mathema\-tic\-al Logic in Paris.
He made a memorisable comment on the  about  the ``wings of logic.''
 Louis Couturat stated that symbolic logic gave wings to mathematics, which involved the
 objection by Henry Poincar\'e that instead of giving it wings logic  had only put mathematics in chains. Robinson rebuffed that however great the mathematician Poncar\'e may had been he was wrong about logic.

The Summer Institute in Logic held in Cornell was one of the most important events for
Abby in~1957, when he had already sought for a job beyond Canada.
Paul Halmos was the initiator of this gathering under the auspices of the American Mathematical Society. He wrote to Edward Hewitt who supervised summer institutes
in the AMS that
logic is a live subject developing rapidly without any support from ``an admiral of the navy
or a~tycoon of industry.'' Leon Henkin and Alfred Tarski backed up Halmos's proposal.
The meeting in Cornell marked the start of the rapid progress of logic in the USA.
Robinson delivered three lectures on relative model completeness and elimination of
quantifiers, on applications of field theory and on proving  theorems ``as done by man, logician, or machine.''
it is curious that Halmos had proclaimed himself to be a~``logician {\it humoris causa}.''
Perhaps, his future invectives against nonstandard analysis demonstrate this status
of his.\footnote{Cp.~\cite[pp.~202--206]{Halmos}.}

In~1957 Robinson had left Canada and returned to his {\it alma mater} in Israel,
where he delivered the compulsory courses on linear algebra and hydrodynamics
and  a special course on logic. In~1959 he was invited to read a course in fluid
mechanics in the Weizmann Institute.
Although Abby's contribution to applied mathematics was fully acknowledged,
the place of his studies in the area had slowly faded out. But
Robinson never lost his admiration of applications. Alec Young, a specialist
in wing theory, remarked that everyone felt that ``the applied mathematician
in Abby was never far away,'' ready to meet any enticing challenge of praxis.
After retirement of Fraenkel, Robinson became the dean of the Mathematics Department
of the Hebrew University of Jerusalem.

In~1960 Robinson spent a sabbatical in Princeton. At the 1960 International Congress
for Logic, Methodology, and Philosophy of Science he made the talk
``Recent Developments in Model Theory'' in which he gave a comprehensive survey
of the pioneering works by Anatoly Maltsev, thus opening them in fact to the logicians
of the USA. Robinson strongly emphasized the importance of the Maltsev  studies demonstrating
how the direct application of model theory led to particular algebraic results.

Soon Robinson was invited to make a plenary talk at the silver anniversary meeting of the Association for
Symbolic Logic which took place on January 24, 1961. This date has become the birthday
of nonstandard analysis.
In summer of 1961  Robinson had been invited to work in the
University of California, Los Angeles
(UCLA), where he moved in July, 1962.
One of the first scholars who shared the ideas of nonstandard analysis was
Wim Luxemburg, an outstanding specialist in functional analysis who primarily
studied Banach lattice theory.\footnote{Cp.~\cite{Lux69}, \cite{Lux72}.}
In may of 1962 Robinson wrote to Luxemburg: ``For some time now I have been thinking
about problems in Functional Analysis but so far as I can see our activities also may
intersect there. Altogether, so far as my standard duties permit, I am now living
in a non-standard mathematical world...\,.'' That was the manner he lived up
to the end of life.

Robinson tried to demonstrate to richness of the new ideas
in most diverse problems. He wrote on the technique of nonstandard analysis
in theoretical physics, studied
nonstandard points on algebraic curves, developed applications of the tools to large
exchange economies, to integration of differential equations, to summations methods,
and so on.\footnote{Cp.~\cite{Rob-2}--\cite{Rob96}.}

In the~1960s Robinson ranked as one of the most popular figures of the
mathematical community. In~196 he was in the center of attention of the participants
of the first international conference on nonstandard analysis which was arranged by Luxemburg in Caltech.
In 1967  Robinson's book~\cite{Rob63} was translated into Russian. But his
{\it magna  opus} on nonstandard analysis was never published in Russian partly in view of
the rise of antisemitism in the academic community of the USSR in those years.

In~1968  Robinson was invited by Nathan Jacobson
to leave UCLA for Yale, where Abby became a tutor of a large group
of young logicians. In~1970. Abby made an invited plenary talk at the International  Congress
of Mathematicians in Nice on ``Forcing in Model Theory.'' In 1971 he received
Sterling Professorship, delivered  a Hedrick lecture at the summer meeting of the
Mathematical Association of America, made a talk at the Fourth International Congress
for Logic, Methodology, and Philosophy of Science in Bucharest,  etc.
In~1972 Abby was elected to the American Academy of Arts and Science,
and in 1973  the  Dutch Mathematical Society decorated Abby with the Brouwer  Medal.

The contributions by Robinson were highly appraised by
the logic genius of the twentieth century Kurt G\"odel who
saw  Robinson as his successor in Princeton.
G\"odel wrote:
``there are good reasons to believe that nonstandard analysis, in
some version or other, will be the analysis of the future'', remarking
that ``his theory of infinitesimals and its application for the solution
of analytical problems seems to me of greatest
importance.''\footnote{Cp.~\cite{Goedel_4-5}.}
Unfortunately, Robinson could not move to Princeton. In November of~ 1973
Abby had begun to fell strong stomach pains, and the doctors found a nonoperable
cancer of the pancreas. Robinson passed away on April 11, 1974 at the age of 55.

\section*{Models of Nonstandard Universes}

Many incomprehensible and obfuscating words can be said
about the nonstandard models of set theory and the
methods of analysis which are based on them.
Some other possibilities are open that we will pursue.

From antiquity mathematics bases on points and numbers.
The most ancient and important method of research consists
in representing numbers by points.
This is the simplest example of modeling for study the properties
of some objects (numbers) by the others---depicting numbers by points.

Let us elaborate this example by recalling
the definitions of points and numbers which are given by Euclid
in his {\it Elements}.
By Definition~I of Book VII
a~monad is ``that by virtue of which
each of the things that exist is called one.''
Euclid proceeds with Definition~II:
``A number is a multitude composed of monads.''
Note that the present-day translations of the Euclid treatise
substitute ``unit'' for ``monad.'' In fact Euclid used the term
{\it M\'{o}}$\nu\alpha\zeta$.

If we look at these definitions from the modern set-theoretic positions
we can assert that two naturals are  a pair of sets one of which includes the other.
But two distinct points has no common elements and so are disjoint as sets.
The points of Euclid are the predecessors of the modern representatives of
the empty set which are called ``atoms.''

So, {\it the presentation of numbers by points on the straight lines
does not preserve the relations of modeled objects as sets}---intersections are not
retained. In other words, the membership between sets is not modeled properly.

A set theoretic model is called  {\it nonstandard},
if the images of modeled objects fail to preserve membership.
This definition belongs to Leon Henkin.
Therefore, the common method of depicting numbers by points
is an example of nonstandard set theoretic modeling.

Let us look at this angle on some classical examples of modeling in mathematics:

(a) The {\it Poincar\'e model} of non-Euclidean geometry
preserves intersections of straight lines and is standard in this respect.

(b) The classical presentation of a {\it separable Hilbert space\/}
as the  $l_2$  sequence space and the $L^2(-\pi,\pi)$ space of square-integrable functions
gave an example of nonstandard modeling. Indeed, the functions $\sin$ and $\cos$
intersect as sets, and the abscissas of the intersection points are interesting.
The presentations of this functions in $l_2$  has an uninteresting intersection.

Therefore, the Riesz--Fisher Isomorphism Theorem for Hilbert spaces
of the same Hilbert dimension uses in fact a nonstandard model, which surely does not
diminish its value and beauty.

Thus example reveals the main particularity of
nonstandard modeling. Nonstandard models are very sensible in regard to what and how
we model and verify claims. In other words, verification, i.~e. discrimination between true and false
propositions, must be explicitly defined.

The above demonstrates that the {\it nonstandard methods of analysis},
i.~e., the techniques based on simultaneous consideration
of standard and nonstandard set-theoretic models, is not something especial, indecent,
or brand-new,
We speak prose but the statement of the fact does make literature critics blush.
We may say that today the {\it possibilities of nonstandard modeling are elaborated
much more than we did it before}.

Return to the most popular and important
nonstandard method of analysis---representation of numbers by points.
Technically speaking, we depict a number~$a$ by the singleton~$\{a\}$.
Clearly, the numbers, say  4 and 5, intersect, whereas their images
$\{4\}$ and $\{5\}$ are disjoint.

Proceed by analogy with all sets, the elements
of von Neumann universe~ ${\mathbb V}$:
$$
\ast:\ a\in{\mathbb V}\rightsquigarrow A\!:=\{a\}\in {\mathbb V}.
$$
As usual,
$$
\mathbb V:=\bigcup\limits_{\alpha\in\On}V_{\alpha},
$$
where
$$
V_{\alpha}\!:=\bigcup\limits_{\beta\in\alpha}\mathscr P(V_{\beta}),
$$
i.~e.,
$V_{\alpha}=\{x \mid (\exists\,\beta )\
(\beta\in\alpha\wedge x\subset V_{\beta})\}$.

\medskip
\centerline{{\epsfxsize8cm\epsfbox{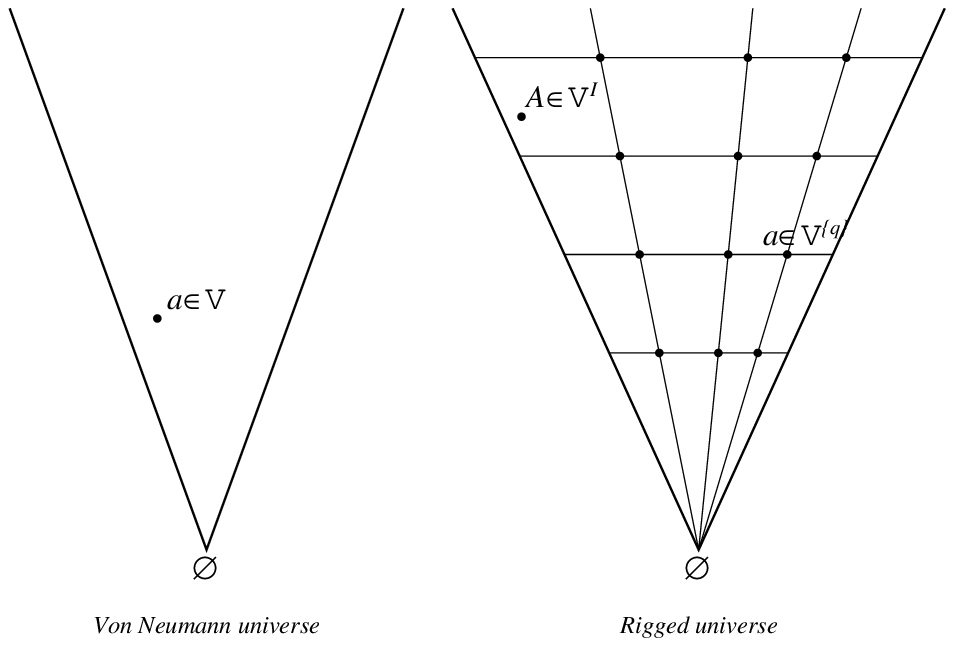}}}
\medskip
\medskip

There is another way of viewing the construction. Consider
$$
{\mathbb V}_q:=\{(q,a)\mid \ a\in{\mathbb V}\}={\mathbb V}^{\{q\}}.
 $$
If $\varphi$  is a proposition about sets, then we have in the model
the {\it transfer principle}:
$\varphi((1,A))\leftrightarrow\varphi(A)$
(i.~e., $(1,a)\subset(1,b)\leftrightarrow a\subset b$ etc.).

Now some abstraction of the construction is in order.
Put
$$
{\mathbb V}_Q:={\mathbb V}^Q:=\{\varphi \mid \varphi :Q\to {\mathbb V},\ \dom\varphi =Q,
\ \im\varphi\subset{\mathbb V}\}.
$$

Let $\mathscr A=\mathscr A(\,\cdot\,)$  be a member of ${\mathbb V}_Q$ and
let $\varphi $ be a formula of set theory.
The law of coordinate-wise modeling is obvious:
$$
\varphi (\mathscr A)\leftrightarrow(\forall\,q\in Q)\ \varphi (\mathscr A(q)).
$$
With this rule ${\mathbb V}_Q$  is clearly a model of set theory.
But it seems that such ``layer-wise'' models brings about nothing new
(despite the fact that the model is nonstandard).

Inspect ${\mathbb V}_Q$ with more attention.
To this end put
$$
[\![\varphi (\mathscr A)]\!]:=\{q\in Q \mid \varphi (\mathscr A(q))\}.
$$
We have the {\it truth value\/} of a formula $\varphi$
as the function $[\![\,\cdot\,]\!]$ with values in~$\mathscr P(Q)$.
Clearly, the {\it transfer principle} holds:
$$
(\varphi \text{ is a theorem})\leftrightarrow[\![\varphi ]\!]=Q.
$$

The following relations are easy:

$$
%\allowdisplaybreaks
%\gather
[\![\varphi\wedge\psi ]\!]=[\![\varphi]\!]\cap[\![\psi]\!],
%\\
$$
$$
[\![\varphi\vee\psi ]\!]=[\![\varphi ]\!]\cup [\![\psi ]\!],
[\![(\forall\,x)\ \varphi]\!]=\bigwedge\nolimits_x
[\![\varphi(x)]\!],
%\\
$$
$$
[\![(\exists\,x)\ \varphi]\!]=\bigvee\nolimits_x[\![\varphi(x)  ]\!],
%\endgather
$$

Moreover, we have the {\it maximum principle}:
$$
(\exists\,x)\ \varphi (x)\to(\exists\,x)\ [\![\varphi (x)]\!]=Q.
$$

Look at $\mathscr R:=\mathbb R^Q$.
Clearly,
$$
[\![\mathscr R\text{ is the field of reals}]\!]=Q;
$$
i.~e., $\mathscr R$ models the reals {\it inside~${\mathbb V}_Q$}.

In other words,
$$
{\mathscr R\!\downarrow}:=\{z\in{\mathbb V}_Q \mid [\![z\in\mathscr R]\!]=Q\}=\mathbb R^Q,
$$
i.~e. The {\it descent} of the reals inside ${\mathbb V}_Q$
is the space of real functions on~$Q$.

We can somehow abstract this construction by replacing~$Q$
with the Stone space of some complete Boolean algebra~$B$,
and $\mathscr P(Q)$, with the clopen algebra of~$Q$.
Strictly speaking, to define the new nonstandard universe~${\mathbb V}^{(B)}$,
put
$$
%\gather
{\mathbb V}_{\alpha}^{(B)}:=\{x \mid  (\exists \beta\in\alpha)\
x:\dom (x)\rightarrow
B\ \wedge\ \dom (x)\subset {\mathbb V}_{\beta}^{(B)}  \},
%\\
$$
$$
{\mathbb V}^{(B)}:=\bigcup\limits_{\alpha\in\On} {\mathbb V}_{\alpha}^{(B)},
%\endgather
$$
where $\alpha$ ranges over the class of all ordinals.

The Boolean truth value $[\![\varphi ]\!]$
is defined by recursion on the complexity of~$\varphi $ while
using the natural interpretation of the logical connectives
and quantifies.
The truth values of the atomic formulas $x\in y$ and $x=y$ for $x,y\in{\mathbb V}^{(B)}$
are defined by transfinite recursion:
$$
%\gather
[\![x\in y]\!]:=\bigvee\limits_{z\in\dom (y)} y(z)\wedge [\![z=x]\!];
%\\
$$
$$
[\![x=y]\!]:=\bigwedge\limits_{z\in\dom (x)} x(z)\Rightarrow [\![z\in y]\!]
\wedge\bigwedge\limits_{z\in\dom (y)} y(z)\Rightarrow [\![z\in x]\!].
%\endgather
$$
(Here
$\Rightarrow$
stands for implication in~$B$.), i.~e. $a\Rightarrow b:=a^\perp\vee b$.)

This we arrive at the Boolean valued universe ${\mathbb V}^{(B)}$
modeling the von Neumann universe~${\mathbb V}$.
The reals~$\mathscr R$ inside the new models has as
its descent the universally complete Kantorovich space
${\mathscr R\!\downarrow}$ whose base in isomorphic to the
initial Boolean algebra~$B$.
This leads to justification of the {\it Kantorovich heuristic principle},
and the theory of Boolean valued models turns into one of the powerful tools
of vector lattice, a classical section of functional analysis.

The theory of Boolean valued models stems from
Dana Scott, Robert Solovay,  Petr Vop\v enka, and Gaisi Takeuti.
These models naturally embrace the nonstandard models by Robinson.

But the importance of infinitesimal methods
is so great that nonstandard analysis takes a rather special place
in mathematics.

Modeling nonstandard analysis, we pass from the usual von Neumann universe
to the ``rigged'' unverse~${\mathbb V}^I$  of the so-called
{\it internal\/} sets with marked frame points,
{\it standard\/} sets, which comprise the copy
${\mathbb V}^S$ of~$\mathbb V$ .
Further analysis shows that ${\mathbb V}^I$ lies in another class, the universe
${\mathbb V}^E$ of {\it external\/} sets satisfying Zermelo axioms.
Some universe ${\mathbb V}^C$
 of~{\it classical\/} sets is distinguished in~ ${\mathbb V}^E$ which is
another realization of the universe of standard sets.
Precisely speaking, there is available some $*$-mapping that
identifies ${\mathbb V}^C$ and ${\mathbb V}^S$ element-wisely.
By analogs of transfer principles, ${\mathbb V}^C$, ${\mathbb V}^S$, and
${\mathbb V}^I$ may be treated as hypostases of the von Neumann universe~${\mathbb V}$.
This complicated interaction of nonstandard models is the background
of modern reconsideration and enrichment of the ancient infinitesimal methods.

%\bigskip
\medskip
\centerline{{\epsfxsize8cm\epsfbox{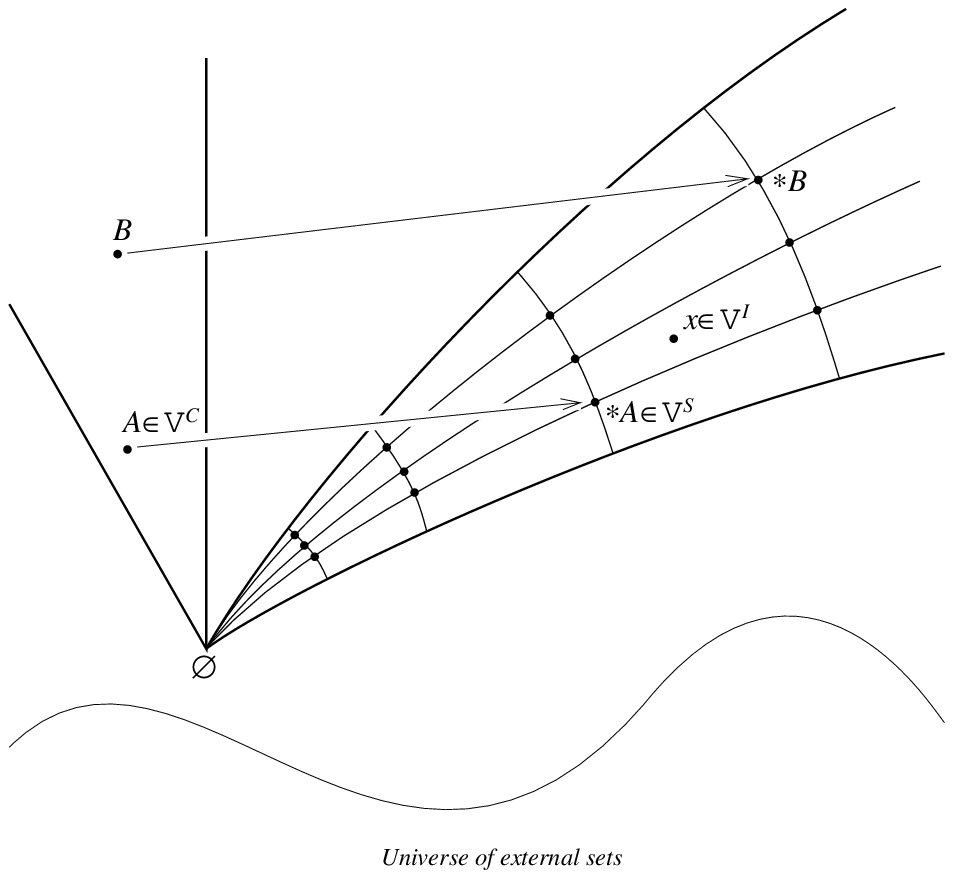}}}
\medskip
%\medskip

\section*{Place of Nonstandard Analysis}

Robinson's nonstandard analysis summarizes the two millennia of the history of
views on actual and potential infinities and paves way to the future of the classical
calculus, suggesting a new paradigm of foundations which is free from
many restrictions of categoricity, prejudice, and bias.

Nonstandard analysis is understood rather broadly today and considered as the branch of
mathematics that bases on nonstandard models of set theory.
In actually this means that under study are two interacting models simultaneously.
Many versions of nonstandard analysis are constructed axiomatically.
The most popular are Nelson's Internal Set Theory and Kawai's Nonstandard
Set Theory.\footnote{Cp.~\cite{Nels77}--\cite{Infa}.}
These theories formalized the ideas
that stem from the ancient views of distinctions between actual and potential infinity.
The theories are conservative extension of Zermelo--Fraenkel Set Theory, having
the same status of justification and rigor in used as for foundations of
mathematics. But the new theories provide new expressible opportunities for
modelling, analyzing, and solving theoretical and applied problems.

The principal starting point of an axiomatic of nonstandard analysis is the
conception that there are two  types of origin objects. The elements of one type are
available to us either immediately or by a potentially infinite processes in the sense
that we can describe such an element directly or prove its unique existence using
already available objects. The available elements are called standard and the
others are called nonstandard. Nonstandard analysis postulates that each infinite
collection of objects has at least one nonstandard element, and every collection of
standard elements is itself standard. This implies the transfer principle that asserts the
cognizability of a standard mathematical property of any standard collection
from the properties of its standards members.

It is worth emphasizing that nonstandard analysis in axiomatic form uses
some new primary concept, the property of every element to be or not to be standard.
In the~``standard''---classical---mathematics of today this property cannot be expressed
and so we cannot speak about actual infinitely large or infinitely small elements.
Moreover, the formal theory of nonstandard analysis is a conservative extension
of the classical set theory. This means that every proposition of classical mathematics
which is proved by using nonstandard analysis can be demonstrated without them new nonstandard tools.
This leads to the popular misconception that nonstandard analysis can add nothing new.
In fact nonstandard analysis is capable of studying the properties
of actual infinites and infinitesimals, suggesting new methods of modeling
or illuminating
the methods of the creators of the calculus like Newton, Leibniz, and Euler who used actual infinities and discriminated
between the assignable and nonassignable reals.

All in all, nonstandard analysis  opens up some new opportunities that are unavailable
in ``standard'' mathematics. In other words, nonstandard analysis
studies the same objects as common mathematics, but it sees in each object
The method of nonstandard analysis is sometimes compared with color TV.
A black and white TV set can see the same objects as a color TV set, but
it cannot discern the variety of colors of the constituents of these elements.
This analogy illustrates the fact that the role of nonstandard analysis
is much broader that providing extra tools for simplifying the apparatus of
standard mathematics by the technique of  ``killing quantifies.''
Nonstandard analysis reveals the rich inner structure of the classical mathematical objects
that are filled with various available and imaginable elements.

To survey the impact of the ideas of nonstandard analysis
is as impossible as to survey applications of  calculus or probability.
Robinson's formalism finds applications in mathematical economics, management,
programming, hydrodynamics, optimization, and elsewhere.
The formalism of Nelson's Internal Set Theory has expanded and enriched
the methodology and the sphere of applications of nonstandard analysis.
The new paradigm is connected with reconsideration of the view of the
classical continuum. In Nelson's theory infinitesimal reside within the unit
segment rather than in ir nonstandard extension as in Robinson's theory.
It is impossible to leave unmarked the rehabilitation of
the Mises frequency approach which was implemented by Edward Nelson in his
conception of ``radically elementary probability theory.''
The external set theories by Toru Kawai and Karel
Hrb\'a\v cek enriched the descriptive and technical  possibilities of nonstandard analysis,
combining the merits of Robinson's and Nelson's formalisms.

Mathematics must constantly fit itself to
the common paradigms of science. Nonstandard analysis
crowns the old-fashioned ideas of the ancient atomism
in much the same way as the Lobachevsky geometry terminated
the dogmatical period of the development of Euclidean geometry.
Robinson suggested a new outlook on the history of mathematical concepts
that underlie the foundations of analysis. His views and approaches
proliferate rapidly these days.\footnote{Cp.~\cite{Notices}.}

The twentieth century is marked with liberation of  humankind
from dogmatism and uniformity.
Filled with the inflammable mixture of genius and villainy of the
population of {\it Homo Sapiens}, the twentieth century will
remain in history the age of liberation of humankind from fascism,
categoricity, absolutism, and domination rather than the age of hatred
and cannibalism. Nonstandard analysis is a produce and
source of freedom.

Humankind will never waste out its intellectual treasures.
Thus there is no doubt that the G\"odel forecast of the future of nonstandard
analysis will turn out prophetic, and some version of nonstandard analysis
will take place of the classical differential and integral calculus of today.
Differentiation as search of trends and integration as prediction of future from trends
are immortal technologies of  mind. New technologies are awaiting humankind
which will use the whole of mathematics in portions incomprehensible today.
This will be the analysis of the future G\"odel had written about.

\end{document}